
\input amstex \documentstyle {amsppt}
\magnification=1200
\hsize=15truecm
\baselineskip=14pt
\hoffset=0truecm

\redefine\cM{{\Cal M}}

\redefine\cF{{\Cal F}}
\redefine\cV{{\Cal V}}
\redefine\cT{{\Cal T}}

\redefine\ee{\Varepsilon}

\redefine\CC{\text{{\bf C}}}
\redefine\NN{\text{{\bf N}}}
\redefine\HH{{\Cal H}}

\redefine\ff{\varphi}

\redefine\HH{{\Cal H}}
\redefine\lb{\lbrack}
\redefine\rb{\rbrack}
\redefine\l{$\lbrack$}
\redefine\r{$\rbrack$}
\redefine\cF{{\Cal F}}
\redefine\cV{{\Cal V}}

\redefine\ee{\epsilon}

\redefine\ct{{\text{{\bf t}}}}
\redefine\SS{\text{S}}
\redefine\SSn{\text{S}_n}
\redefine\id{\text{{\bf 1}}}
\redefine\subnot{\subset}
\redefine\cFT{\cF_T}
\redefine\la{\langle}
\redefine\ra{\rangle}
\redefine\raT{\rangle_T}
\redefine\Pn{P^{(n)}}
\redefine\Pnpe{P^{(n+1)}}
\redefine\Rnpe{R^{(n+1)}}
\redefine\Hs{\tilde\HH}
\redefine\HsFT{\Hs\otimes\cF_T}
\redefine\dim{\text{dim}}
\redefine\qii{q_{ii}}
\topmatter
\title COMPLETELY POSITIVE MAPS ON COXETER GROUPS,\\
DEFORMED COMMUTATION RELATIONS,\\ AND OPERATOR SPACES \endtitle
\author {\bf Marek Bo\D zejko$^1$ and Roland Speicher$^2$}\\ \quad\\
$^1$ Instytut Matematyczny\\
Uniwersytet Wroc\l awski\\ Plac Grunwaldzki 2/4\\
50-384 Wroc\l aw\\ Poland \\ e-mail: bozejko \@ plwruw11.bitnet\\
\quad\\ and\\ \quad\\
$^2$ Institut f\"ur Angewandte Mathematik\\ Universit\"at Heidelberg\\
Im Neuenheimer Feld 294\\ D-69120 Heidelberg\\ Federal Republic of Germany\\
e-mail: L95 \@ vm.urz.uni-heidelberg.de 
\endauthor
\endtopmatter
\document
\heading
{\bf Abstract}
\endheading
In this article we prove that quasi-multiplicative (with respect to
the usual length function)
mappings on the permutation group $\SSn$ (or, more generally, on arbitrary
amenable Coxeter groups), determined by self-adjoint contractions fulfilling
the braid or Yang-Baxter relations,
are completely positive. We point out the
connection of this result with the construction of a Fock
representation of the deformed commutation relations
 $d_id_j^*-\sum_{r,s} t_{js}^{ir} d_r^*d_s=\delta_{ij}\id$,
where the matrix $t_{js}^{ir}$ is given by a self-adjoint contraction
fulfilling the braid relation.
Such deformed commutation relations give examples for
operator spaces as considered by Effros, Ruan and Pisier. The corresponding
von Neumann algebras, generated by $G_i=d_i+d_i^*$, are typically not
injective.
\par
\pagebreak
\baselineskip=12pt
\heading
{\bf 1. Introduction}
\endheading
We will prove in this paper the following result.
\proclaim{Theorem 1.1} Consider for fixed $n\in\NN$ the permutation
group $\SSn$ and denote by $\pi_i\in\SSn$ ($i=1,\dots,n-1$) the
transposition between $i$ and $i+1$. Furthermore, let operators
$T_i\in B(\HH)$ ($i=1,\dots,n-1$) on some Hilbert space $\HH$ be given
with the properties:\newline
i) $T_i^*=T_i$ for all $i=1,\dots,n-1$\newline
ii) $\Vert T_i\Vert\leq 1$ for all $i=1,\dots,n-1$\newline
iii) The $T_i$ satisfy the braid relations:
 $$T_iT_{i+1}T_i=T_{i+1}T_iT_{i+1}\qquad\text{for all $i=1,\dots,n-2$}$$
 $$T_iT_j=T_jT_i\qquad\text{for all $i,j=1,\dots,n-1$ with $\vert i-j\vert
\geq 2$}$$
Define now a function
 $$\ff:\SSn\to B(\HH)$$
by quasi-multiplicative extension of
 $$\ff(e)=\id,\qquad \ff(\pi_i)=T_i,$$
i.e. for a reduced word $\SSn\owns \sigma=\pi_{i(1)}\dots\pi_{i(k)}$
we put $\ff(\sigma)=T_{i(1)}\dots T_{i(k)}$.\newline
Then $\ff$ is a completely positive map,
i.e. for all $l\in\NN$, $f_i\in\CC\SSn$, $x_i\in\HH$ ($i=1,\dots,l$) we
have
 $$\langle \sum_{i,j=1}^l\ff(f_j^*f_i)x_i,x_j\rangle\geq 0.$$
\endproclaim
One should note that the braid relations of the $T_i$
ensure \l Bou\r\ that $\ff$ is well defined. Of course, we do not
assume that $T_i^2=\id$. In this case, $\ff$ were a representation of
$\SSn$ and the theorem would be trivial.\par
To get a flavour of the meaning of this theorem let us just mention, that
for $\SS_3$ the statement is equivalent to the fact that the operator
 $P=\id+T_1+T_2+T_1T_2+T_2T_1+T_1T_2T_1$
is strictly positive, whenever $T_i^*=T_i$, $\Vert T_i\Vert<1$, and
$T_1T_2T_1=T_2T_1T_2$.
\par
Theorem 1.1 is also valid much more generally, namely for all finite
(or even amenable) Coxeter groups. The formulation and the proof of this
generalization will be given in Sect. 1.
\par
The motivation for our Theorem 1.1 comes from investigations on
perturbed commutation relations. The crucial step in establishing the
existence of a Fock representation of such relations is the
positivity of some map on the permutation groups $\SSn$. So, in
\l BSp1\r\, we investigated the relations
 $$c_ic_j^*-q c_j^*c_i=\delta_{ij}\id$$
for a real $q$ with $\vert q\vert\leq1$, and we needed essentially the fact
that
 $$\ff:\SSn\to\CC,\qquad \pi\mapsto q^{\vert\pi\vert}$$
is a positive definite function for all $n$, where $\vert\pi\vert$
denotes the number of inversions of $\pi$. Different proofs of this
positive  definiteness, and correspondingly of the existence of the
Fock representation of the $q$-relations, are now available, see
\l BJS,BSp1,BSp2,Spe2,Fiv,Gre,Zag\r.
\par
In \l Spe2\r, we considered, more generally, the relations
 $$d_id_j^*-q_{ij}d_j^*d_i=\delta_{ij}\id$$
for $-1\leq q_{ij}=q_{ji}\leq 1$ and proved by central limit arguments the
existence of a Fock representation. In Sect. 3, we will construct
the Fock representation of these deformed commutation relations,
now even for the most general case of complex $q_{ij}$ with $\bar q_{ij}=
q_{ji}$. We will see that again the positivity of some map on $\SSn$
is the key point behind this construction. This positivity will then follow
as a special case of our general Theorem 1.1.\par
Our construction of the $q_{ij}$-relations depends essentially on
some operator $T$, which is a self-adjoint contraction and fulfills the
braid or Yang-Baxter relation. Thus, our natural frame in Sects. 3 and 4
will be that we consider the general deformed commutation relations
 $$d_id_j^*-\sum_{r,s} t_{js}^{ir} d_r^*d_s=\delta_{ij}.$$
Such general Wick ordering relations are
also investigated by
J\o rgensen, Schmitt, and Werner \l JSW2\r. Whereas in the most general
case, without any assumptions on $t_{js}^{ir}$
apart from the necessary $\bar t_{ab}^{dc}=
t_{dc}^{ab}$, nothing can be said about the existence of a Fock
representation, we get, by Theorem 1.1, a proof for the existence of this
representation in the case
where the matrix $t_{js}^{ir}$ is given by a self-adjoint contraction $T$
fulfilling the braid relation.
\par
In Sect. 4, we examine the deformed commutation relations from an operator
space point of view, namely we extend a result of Haagerup and Pisier
and show that the operator space generated by the $G_i:=d_i+d_i^*$
is completely isomorphic to the canonic operator space $R\cap C$, which
means
 $$
   \Vert(a_1,\dots,a_N)\Vert_{\max}
\leq\Vert\sum_{i=1}^N a_i\otimes G_i\Vert\leq \frac 2{\sqrt{1-q}}
   \Vert(a_1,\dots,a_N)\Vert_{\max}$$
for all bounded operators $a_1,\dots,a_N$ on some Hilbert space, where
 $$\Vert(a_1,\dots,a_N)\Vert_{\max}:=
   \max(\Vert\sum_{i=1}^N a_ia_i^*\Vert^{1/2},
\Vert\sum_{i=1}^N a_i^*a_i\Vert^{1/2}).$$
We will also make some remarks on the von Neumann algebra generated by
the $G_i$. In particular, we show that it is typically not injective.
\par
Our main theorem, 1.1 and its general version 2.1, considers operator
valued functions which are quasi-multiplicative with respect to the
usual length function (=minimal number of generators). In Sect. 5, we
replace this length function by another, also quite natural one (=
minimal number of {\it different} generators) and prove the analogue of
2.1 for this case.
\newline
\vskip1cm
\heading
{\bf 2. Completely positive maps on finite Coxeter groups}
\endheading
Let $(W,S)$ be a Coxeter system consisting of a Coxeter group $W$ and
a set $S=\{s_1,\dots,s_n\}$ of generators. This means that $W$ is the
group generated by the elements $s_i=s_i^{-1}\in S$ and that for each
two distinct generators $s_i,s_j\in S$ ($i\not= j$) there exists a
natural number $m_{ij}\geq 2$ such that we have the relation
 $$(s_is_j)^{m_{ij}}=e,$$
where $e$ is the unit element of $W$. The fact $s_i=s_i^{-1}$ can also
be stated in this form as $m_{ii}=1$. In the following we will only consider
finite Coxeter groups $W$.\par
For each $\sigma\in W$ we denote by $\vert \sigma\vert$ the length of
$\sigma$ with respect to $S$, i.e.
 $$\vert\sigma\vert:=\min\{k\in\NN\mid\text{there exist $s_{i(1)},\dots
,s_{i(k)}\in S$ with $\sigma=s_{i(1)}\dots s_{i(k)}$}\},$$
and $\vert e\vert =0$.\par
The example of $\SSn$ fits into this frame by putting
 $W=\SSn$, $S=\{\pi_1,\dots,\pi_{n-1}\}$. The length function $\vert\pi\vert$
is then given by the number of inversions and
the relations are
given by $m_{ij}=2$ for $\vert i-j\vert\geq2$, i.e.
 $$\pi_i\pi_j\pi_i\pi_j=e \Leftrightarrow \pi_i\pi_j=\pi_j\pi_i\qquad
(\vert i-j\vert\geq 2)$$
and $m_{i,i+1}=3$, i.e.
 $$\pi_i\pi_{i+1}\pi_i\pi_{i+1}\pi_i\pi_{i+1}=e \Leftrightarrow
\pi_i\pi_{i+1}\pi_i=\pi_{i+1}\pi_i\pi_{i+1}.$$
For a general Coxeter group $W$, we will also rewrite the defining
relations $(s_is_j)^{m_{ij}}$ $=e$ in the braid like form
 $$\undersetbrace \text{$m_{ij}$ factors} \to {s_is_js_i\dots}=
   \undersetbrace \text{$m_{ij}$ factors} \to {s_js_is_j\dots},$$
which means
 $$(s_is_j)^{m_{ij}/2}=(s_js_i)^{m_{ij}/2}\qquad\text{for $m_{ij}$ even}$$
and
 $$(s_is_j)^{(m_{ij}-1)/2}s_i
=(s_js_i)^{(m_{ij}-1)/2}s_j\qquad\text{for $m_{ij}$ odd.}$$
Let now self-adjoint contractions $T_i\in B(\HH)$
on some Hilbert space be given which fulfill the generalized braid
relations
 $$\undersetbrace\text{$m_{ij}$ factors}\to{T_iT_jT_i\dots}=
   \undersetbrace\text{$m_{ij}$ factors}\to{T_jT_iT_j\dots}$$
for all $i,j=1,\dots,n$ with $i\not= j$.
Then we define the mapping
 $$\ff:W\to B(\HH)$$
by
 $\ff(e)=\id$ and
 $$\ff(\sigma)=T_{i(1)}\dots T_{i(k)}\quad\text{for}\quad
\sigma=s_{i(1)}\dots s_{i(k)} \quad\text{with $\vert \sigma\vert=k$.}$$
It is known \l Bou\r\ that the generalized braid relations for the $T_i$
ensure that this definition of $\ff$ is well defined.
We can also state our definition in the way that we put $\ff(s_i)=T_i$ and
extend $\ff$ in a quasi-multiplicative way, which means
 $$\ff(\sigma_1\sigma_2)=\ff(\sigma_1)\ff(\sigma_2)\qquad\text{if}\qquad
\vert\sigma_1\sigma_2\vert=\vert\sigma_1\vert+\vert\sigma_2\vert.$$
Note that the self-adjointness of the $T_i$ implies
 $\ff(\sigma^{-1})=\ff(\sigma)^*$.\par
Let us extend $\ff$ from $W$ to its group algebra
 $$\CC W:=\{f=\sum_{\sigma\in W}f(\sigma)\delta_\sigma\}$$
(with the usual multiplication
($\delta_\sigma\delta_\pi=\delta_{\sigma\pi}$)
and involution ($\delta_\sigma^*=\delta_{\sigma^{-1}}$) structure) in
the canonical way
 $$\ff(\sum_{\sigma\in W} f(\sigma)\delta_{\sigma})=\sum_{\sigma\in W}
f(\sigma)\ff(\sigma),$$
then we can state our main result in the following way.
\proclaim{Theorem 2.1} Let $T_i\in B(\HH)$ ($i=1,\dots,n$) be bounded operators
on some Hilbert space $\HH$, which fulfill the following assumptions:
\roster
\item"i)"
$T_i^*=T_i$ for all $i=1,\dots,n$
\item"ii)"
$\Vert T_i\Vert\leq 1 $ for all $i=1,\dots,n$
\item"iii)"
We have for all $i,j=1,\dots,n$ with $i\not=  j$ the generalized braid
relations
 $$\undersetbrace \text{$m_{ij}$ factors}\to{T_iT_jT_i\dots}=
   \undersetbrace \text{$m_{ij}$ factors}\to{T_jT_iT_j\dots}.$$
\endroster
Then the quasi-multiplicative map
 $$\ff:\CC W\to B(\HH)$$
given by
 $$\ff(e)=\id,\qquad \ff(s_i)=T_i\quad (i=1,\dots,n)$$
is completely positive, i.e. for all $l\in\NN$, $f_i\in\CC W$, $x_i\in\HH$
($i=1,\dots,l$) we have
 $$\langle \sum_{i,j=1}^l \ff(f_j^*f_i)x_i,x_j\rangle\geq 0.$$
\endproclaim
\demo{Remark} Another equivalent characterization of complete positivity
is the fol\-low\-ing (see, e.g., \l Pau\r): For arbitrary $\alpha:W\to\HH$
(with finite support) we have
 $$\sum_{\rho,\sigma\in W} \la\ff(\rho^{-1}\sigma)\alpha(\sigma),\alpha
(\rho)\ra\geq 0.$$
This formulation is the operator valued version of the definition of a
positive definite function on $W$.
\enddemo
Theorem 2.1 follows essentially from the following theorem.
\proclaim{Theorem 2.2}
Let $T_i$ and $\ff$ be as in 2.1. Then the operator
 $$P:=\sum_{\sigma\in W}\ff(\sigma)$$
is positive, i.e. $P\geq 0$.
\endproclaim
By putting all $x_i\equiv x$ it is clear that complete positivity implies
$P\geq 0$. Let us now see how we get, in the other direction, 2.1 from
2.2.
\demo{Proof of 2.2 $\Rightarrow$ 2.1}
Let $\lambda$ be the left regular representation of $W$ acting on
$l^2(W)\hat =\CC W$ equipped with the scalar product
 $\langle f,g\rangle=\sum_{\sigma\in W} \bar f(\sigma) g(\sigma)$,
i.e.
 $$(\lambda(\rho)f)(\sigma):=f(\rho^{-1}\sigma)\quad\text{or}\quad
\lambda(\rho)\delta_\sigma=\delta_{\rho\sigma}
\qquad\text{for $\rho,\sigma
\in W$, $f\in l^2(W)$.}$$
If we now define the operators
 $$\hat T_i:=\lambda(s_i)\otimes T_i\qquad\text{on $l^2(W)\otimes\HH$},$$
then they also satisfy the assumptions of Theorem 2.2, which yields
 $$\hat P:=\sum_{\sigma\in W}\hat\ff(\sigma)\geq 0,$$
where $\hat\ff$ is the quasi-multiplicative function given by the $\hat T_i$,
clearly
 $$\hat\ff(\sigma)=\lambda(\sigma)\otimes\ff(\sigma)\qquad\text{for $\sigma\in
 W$.}$$
The positivity of $\hat P$ implies now
 $$\langle \sum_{i,j=1}^l \ff(f_j^*f_i)x_i,x_j\rangle\geq 0$$
in the following way: Put
 $$\alpha:
=\sum_{\rho\in W}\delta_\rho\otimes\bigl(\sum_i f_i(\rho^{-1})x_i\bigr)\in
l^2(W)\otimes\HH.$$
Then
 $$\align
0&\leq \langle \hat P\alpha,\alpha\rangle\\
&=\sum_\sigma\langle \bigl(\lambda(\sigma)\otimes\ff(\sigma)\bigr)
\alpha,\alpha
\rangle\\
&=\sum_{\sigma,\rho,\tau,i,j}\langle \bigl(\lambda(\sigma)\otimes\ff(\sigma)
\bigr)
\bigl(\delta_\rho\otimes f_i(\rho^{-1})x_i\bigr),\delta_\tau\otimes
f_j(\tau^{-1})x_j\rangle\\
&=\sum_{\sigma,\rho,\tau,i,j} \langle \delta_{\sigma\rho}\otimes
f_i(\rho^{-1})\ff(\sigma)x_i,\delta_\tau\otimes f_j(\tau^{-1})x_j\rangle\\
&=\sum_{\sigma,\rho,\tau,i,j}\langle \delta_{\sigma\rho},\delta_\tau\rangle
\langle \bar f_j(\tau^{-1})f_i(\rho^{-1})\ff(\sigma)x_i,x_j\rangle\\
&=\langle \sum_{i,j} \ff(f_j^*f_i)x_i,x_j\rangle.\endalign$$
The last line follows from
 $$\langle \delta_{\sigma\rho},\delta_\tau\rangle=\cases 1,& \sigma=\tau
\rho^{-1}\\ 0,& \text{else}\endcases$$
and the fact that with
 $$f_i=\sum_\rho f_i(\rho^{-1})\delta_{\rho^{-1}},\qquad
   f_j^*=\sum_\tau \bar f_j(\tau^{-1})\delta_{\tau}$$
we have
 $$f_j^*f_i=\sum_{\rho,\tau}\bar f_j(\tau^{-1})f_i(\rho^{-1})\delta_{\tau
\rho^{-1}}.$$
\line{\hfill $\diamondsuit$}
\enddemo
So we are left with the proof of 2.2. Note first that it suffices to treat
the case of strict contractions.
\proclaim{Theorem 2.3}
Let $T_i$ and $\ff$ be as in 2.1, but with the stronger assumption of
strict contractivity, i.e. $\Vert T_i\Vert< 1$ for all $i=1,\dots,n$.
Then the operator
 $$P:=\sum_{\sigma\in W}\ff(\sigma)$$
is strictly positive, i.e. $P>0$.
\endproclaim
Theorem 2.2 can be infered from this version in the following way.
\demo{Proof of 2.3 $\Rightarrow$ 2.2}
Consider $T_i^{(t)}:=tT_i$ ($i=1,\dots,n$) for $0\leq t<1$. If the $T_i$
fulfill the assumptions of 2.2, then the $T_i^{(t)}$ fulfill the assumptions
of 2.3. Thus
 $$P^{(t)}=\sum_{\sigma\in W} \ff(\sigma) t^{\vert\sigma\vert}>0\qquad
\text{for all $0\leq t<1$.}$$
If now $t\nearrow 1$, then $P^{(t)}\to P$ uniformly and we get the assertion.
\newline
\line{\hfill $\diamondsuit$}
\enddemo
To prove 2.3 we reduce the assertion about strict positivity to one
about invertibility. Note that $P$ is self-adjoint, since
 $$P^*=\sum_{\sigma\in W}\ff(\sigma)^*=\sum_{\sigma\in W}\ff(\sigma^{-1})=P.$$
\proclaim{Theorem 2.4}
 Let $T_i$, $\ff$, and $P$ be as in 2.3. Then $P$ is invertible.
\endproclaim
For the reduction of 2.3 to 2.4 we need a fact on the norm-continuity of
the smallest element in the spectrum of a self-adjoint operator. Define
for a self-adjoint operator $A\in B(\HH)$ the number
 $$m_0(A):=\inf\{\la Ax,x\ra\mid x\in\HH, \Vert x\Vert=1\}$$
as the smallest element in the (compact!) spectrum of $A$.
\proclaim{Lemma 2.5}
We have for arbitrary self-adjoint operators $A,B\in B(\HH)$
 $$\vert m_0(A)-m_0(B)\vert\leq \Vert A-B\Vert.$$
\endproclaim
\demo{Proof of 2.5}
Assume $m_0(A)\geq m_0(B)$. Fix an arbitrary $\ee>0$. Then there exists
$x\in\HH$ with $\Vert x\Vert=1$ such that $m_0(B)\geq \la Bx,x\ra-\ee$.
Since $m_0(A)\leq \la Ax,x\ra$ we have
 $$\vert m_0(A)-m_0(B)\vert\leq\la Ax,x\ra-\la Bx,x\ra+\ee=
\la(A-B)x,x\ra+\ee\leq\Vert A-B\Vert+\ee.$$
For $\ee\to 0$ we get the assertion.\newline
\line{\hfill $\diamondsuit$}
\enddemo
\demo{Proof of 2.4 $\Rightarrow$ 2.3}
Consider again the collection of $T_i^{(t)}:=t T_i$ ($i=1,\dots,n$)
for all $t$ with $0\leq t\leq 1$.
Then, by 2.4,
$P^{(t)}=\sum\ff(\sigma)t^{\vert\sigma\vert}$ is invertible for all
$0\leq t\leq 1$ and we have $P^{(0)}=\id$ and $P^{(1)}=P$. Furthermore,
the map $t\mapsto P^{(t)}$ is norm-continuous.
Put $m_0(t):=m_0(P^{(t)})$ (note $P^{(t)*}=P^{(t)}$).
Since
 $$\vert m_0(t_1)-m_0(t_2)\vert\leq \Vert P^{(t_1)}-P^{(t_2)}\Vert,$$
the mapping $t\mapsto m_0(t)$ is continuous. But now invertibility of
$P^{(t)}$ implies $m_0(t)\not=0$. Because of $m_0(0)=m_0(\id)=1$, we have
$m_0(t)>0$ for all $0\leq t\leq 1$, in particular $m_0(P)=m_0(1)>0$, i.e.
$P>0$.\newline
\line{\hfill $\diamondsuit$}
\enddemo
Up to now we have only used very general arguments for the reduction of
our theorem. This reduction has led us to a statement on invertibility of
some operator $P\in\CC W$. This is now an algebraic problem which can be
\lq calculated' in our group algebra. Of course, now we need the special
structure of Coxeter groups. The proof will be by induction on the
cardinality of Coxeter generators of
parabolic subgroups of $W$. \par
For $J\subseteq S$,
let $W_J$ be the subgroup of $W$ generated by all $s\in J$.
Such subgroups are called parabolic. They are also Coxeter groups, given by
the system $(W_J,J)$. We need now the following known facts on these
subgroups (see, e.g., \l Bou,Car\r):
For $J\subseteq S$ we define
 $$D_J:=\{\sigma\in W\mid \vert \sigma s\vert=\vert\sigma\vert+1 \quad
\text{for all $s\in J$}\},$$
i.e. $\sigma\in D_J$ if and only if $\sigma$ is the element of smallest
length in
the coset $\sigma\cdot W_J$. Thus $D_J$ is a canonical representative
system of the cosets of $W_J$. If we define for $\sigma\in W$ the set
 $$J_{\sigma}:=\{s \in S\mid \vert\sigma s\vert=\vert\sigma\vert+1\},$$
then the definition of $D_J$ can also be put in the way
 $$\sigma\in D_J\Longleftrightarrow  J\subseteq J_{\sigma}.$$
This characterization gives at once the Euler-Solomon-formula \l Sol\r\
for all
$\sigma\in W$
 $$\sum\Sb J\subseteq S\\\text{with}\\\sigma\in D_J\endSb (-1)^{\vert J\vert}
=\sum_{J\subseteq J_{\sigma}}(-1)^{\vert J\vert}
=\cases 0,&\text{if $\sigma\not=\sigma_0$}\\
1,&\text{if $\sigma=\sigma_0$},\endcases$$
where $\sigma_0$ is the unique element in $W$ with the greatest length, i.e.
the unique element with the property $J_{\sigma_0}=\emptyset$.\par
Furthermore, we
have the nice property that
each element $\sigma\in W$ can uniquely be written
in the form $\sigma=\tau_J\sigma_J$ with $\tau_J\in D_J$ and
$\sigma_J\in W_J$, and with $\vert\sigma\vert=\vert\tau_J\vert+
\vert\sigma_J\vert$.
\par
Note in particular that, for $J=S$, we have $W_S=W$ and hence $D_S=\{e\}$.
In the next section we will need the following special case: For $W=\SS_{n+1}$
consider $J:=\{\pi_2,\pi_3,\dots,\pi_n\}$, i.e. $W_J\cong S_n$. Then
$D_J=\{e,\pi_1,\pi_2\pi_1,\dots,\pi_{n}\pi_{n-1}\dots\pi_1\}$. As an
example, take $n=2$, then $W=\SS_3=\{e,\pi_1,\pi_2,\pi_1\pi_2,\pi_2\pi_1,
\pi_1\pi_2\pi_1=\pi_2\pi_1\pi_2\}$, $J=\{\pi_2\}$,
$W_J=\{e,\pi_2\}\cong \SS_2$, $D_J=\{e,\pi_1,\pi_2\pi_1\}$.\par
Now define for an arbitrary subset $A\subseteq W$ the operator
 $$P(A):=\sum_{\sigma\in A}\ff(\sigma).$$
Then the uniqueness of the decomposition $W=D_JW_J$ and the
quasi-multiplicity of $\ff$ give for all $J\subseteq S$
 $$P=P(W)=P(D_J)P(W_J).$$
For the above example of $W=\SS_3$ and $J=\{\pi_2\}$ this decomposition
is given by
 $$\id+T_1+T_2+T_1T_2+T_2T_1+T_1T_2T_1=(\id+T_1+T_2T_1)(\id+T_2).$$
The crucial point for our induction is now the translation of the
Euler-Solomon-formula to our operators $P(A)$.
\proclaim{Lemma 2.6}
Let $(W,S)$ be an arbitrary finite Coxeter group and $\sigma_0$ the
unique longest element in $W$. Then we have
 $$\sum_{J\subseteq S} (-1)^{\vert J\vert} P(D_J)=\ff(\sigma_0).$$
\endproclaim
\demo{Proof of 2.6}
We have
 $$\align
 \sum_{J\subseteq S}(-1)^{\vert J\vert} P(D_J)&=\sum_{J\subseteq S}
(-1)^{\vert J\vert}\sum_{\sigma\in D_J}\ff(\sigma)\\
&=\sum_{\sigma\in W}\bigl(\sum\Sb J\subseteq J_{\sigma}
\endSb (-1)^{\vert J\vert}\bigr) \ff(\sigma)\\
&=\ff(\sigma_0).\endalign$$\nopagebreak
\line{\hfill $\diamondsuit$}
\enddemo
\demo{Proof of 2.4}
We prove this by induction on the cardinality of $S$.
If $\vert S\vert=0$, then $W=\{e\}$ and $P=\id$ is invertible.
If $\vert S\vert=1$, then $W=\{e,s_1\}$ and $P=\id+T_1$ is invertible
because of $\Vert T_1\Vert<1$.\newline
Assume now we know the invertibility of $P(\tilde W)$ for all finite Coxeter
groups $(\tilde W,\tilde S)$ with $\vert \tilde S\vert\leq n-1$.
Consider an arbitrary finite
Coxeter group $(W,S)$ with $\vert S\vert=n$.
By $J\subset S$ we will denote in the following the situation
that $J\subseteq S$, but $J\not= S$.
Then we have by induction hypothesis the invertibility of $P(W_J)$ for
all $J\subnot S$, hence
 $$P(D_J)=P(W)P(W_J)^{-1}.$$
Lemma 2.6 yields then
 $$\align
P(W)\bigl\{\sum_{J\subnot S} (-1)^{\vert J\vert} P(W_J)^{-1}\bigr\}&=
\sum_{J\subnot S} (-1)^{\vert J\vert} P(D_J)\\
&=\ff(\sigma_0)-(-1)^{\vert S\vert}P(D_S)\\
&=\ff(\sigma_0)-(-1)^{\vert S\vert}\id.\endalign$$
Since $\Vert \ff(\sigma_0)\Vert<1$, the element
 $\ff(\sigma_0)-(-1)^{\vert S\vert}\id$ is invertible and we get
 $$P(W)\bigl\{\sum_{J\subnot S}(-1)^{\vert J\vert} P(W_J)^{-1}\bigr\}
\bigl\{\ff(\sigma_0)-(-1)^{\vert S\vert}\id\bigr\}^{-1}=\id,$$
hence $P=P(W)$ is right invertible. Because of $P^*=P$ it is also left
invertible, hence invertible.\newline
\line{\hfill $\diamondsuit$}
\enddemo
\demo{Remarks}
1) If we specialize to $W=\SSn$ and $S=\{\pi_1,\dots,\pi_{n-1}\}$ then
we recover Theorem 1.1 from the introduction. Note that even in this
case our main step, namely the positivity of $P$, is by no means trivial.
E.g., for $\SS_3$ it states that the operator
 $P=\id+T_1+T_2+T_1T_2+T_2T_1+T_1T_2T_1$
is positive, whenever $T_i^*=T_i$, $\Vert T_i\Vert<1$, and $T_1T_2T_1=T_2
T_1T_2$.
\newline
2) Theorem 2.1 is also true for amenable Coxeter
groups. Since we know by a result of de la Harpe \l deH\r\ that amenable
Coxeter groups are either finite or
affine Coxeter groups and hence the cardinality of the
set $\{\sigma\mid\vert \sigma\vert\leq k\}$ is
at most of polynomial growth in $k$ (see \l Bou\r\ for the structure of
affine Coxeter groups),
the operator $P$ is also well defined in the case $\Vert T_i\Vert
<1$ for amenable Coxeter groups. In this case, all our arguments remain the
same, only in Lemma 2.6 the value $\ff(\sigma_0)$ on the right side of the
equation has to be replaced by 0 if the Coxeter group is infinite.
Thus we get in the same manner as for
finite groups the assertion of 2.1 also for amenable groups in the case
$\Vert T_i\Vert<1$. Since the statement on complete positivity involves
only finite sums, we can now carry out the limit $\Vert T_i\Vert\nearrow1$
and obtain in this way the validity of Theorem 2.1 for all amenable Coxeter
groups.\newline
3) It is an open question whether 2.1 is true for all infinite Coxeter
groups. What can be proved in this general case is the validity of 2.1 for
all Coxeter groups in the special case of scalar valued $T_i\in\CC$. This
proof uses other methods and will be published elsewhere \l Boz2\r. In the
special case, when $W$ is the free product of 2-elements
groups, i.e. when we
have as only relations $s_i^2=e$ for all $i$, then Theorem 2.1 was also
proved for the general operator valued case, see \l Boz1\r.
\enddemo
\vskip1cm
\heading
{\bf 3. Fock Representation of Deformed Commutation Relations}
\endheading
We will now use our general result 1.1 for the construction of the
Fock representation of the $q_{ij}$-relations
 $$d_id_j^*-q_{ij}d_j^*d_i=\delta_{ij}\id$$
for $\bar q_{ij}=q_{ji}$ and $\vert q_{ij}\vert\leq 1$, i.e. we are looking
for operators $d_i$, $d_i^*$ on some Hilbert space $\tilde\HH$ and some
vector $\Omega\in\tilde\HH$ (called vacuum), such that $d_i$ and $d_i^*$
are adjoints of each other and all annihilation operators $d_i$ annihilate
the vacuum: $d_i\Omega=0$. One should note that these requirements determine
the structure of the Fock representation up to unitary equivalence, the
only problem is to prove the existence of such a structure, in particular
to show the positivity of the corresponding scalar product in $\tilde \HH$.
\par
We will treat in the following a more general case and specify this in the
end to the above mentioned relations. Assume we are given some operator
$T$ and a Hilbert space $\HH$ such that $T\in B(\HH\otimes\HH)$ is a
self-adjoint contraction ($T^*=T$, $\Vert T\Vert\leq 1$) and such that
it fulfills the braid relation
 $$(\id\otimes T)(T\otimes \id)(\id\otimes T)
=(T\otimes\id)(\id\otimes T)(T\otimes\id),\tag\text{BR}$$
where $\id\otimes T$
and $T\otimes\id$ are the natural amplifications of $T$ to
$\HH\otimes\HH\otimes\HH$. Then we define
 $$T_i:=\undersetbrace\text{$i-1$ times}\to {\id\otimes\dots\otimes\id}
\otimes T\qquad\text{acting on $\HH^{\otimes(i+1)}$},$$
and by amplification
also on all $\HH^{\otimes n}$ with $n\geq i+1$. The $T_i$
fulfill the assumptions of Theorem 1.1.
\par
The braid relation (BR) appears also in a lot of contextes under the
name \lq Yang-Baxter equation', see, e.g., \l Man,Jim,Wen\r.\par
Now we define, for each $f\in\HH$, a creation operator $d^*(f)$ and an
annihilation operator $d(f)$ on a dense subset $\cF$ of the full Fock
space
 $\bigoplus_{n=0}^\infty \HH^{\otimes n}$, where $\HH^0:=\CC\Omega$
($\Vert\Omega\Vert=1$), $\cF$ being the set of finite linear combinations of
product vectors. On the full Fock space we have the canonic free
creation and annihilation operators given by (see \l Eva,Voi,Spe1\r)
 $$\align
l^*(f)\Omega&=f\\
l^*(f)f_1\otimes\dots\otimes f_n&=f\otimes f_1\otimes\dots\otimes f_n
\endalign$$
and
 $$\align
l(f)\Omega&=0\\
l(f)f_1\otimes\dots \otimes f_n&=\langle f,f_1\rangle f_2\otimes\dots\otimes
f_n.\endalign$$
We define now our deformation by
 $$d^*(f):=l^*(f)$$
and
 $$d(f):=l(f)(\id+T_1+T_1T_2+\dots+T_1T_2\dots T_{n-1})\qquad
\text{on $\HH^{\otimes n}$.}$$
Of course, $d(f)$ and $d^*(f)$ are not adjoints of each other with
respect to the usual scalar product $\langle\enskip,\enskip\rangle$. Thus
we introduce a new one $\langle\enskip,\enskip\rangle_T$ given by
 $$\langle \xi,\eta\rangle_T:=\delta_{nm}\langle \xi,P^{(n)}\eta\rangle
\qquad\text{for $\xi\in\HH^{\otimes n}$, $\eta\in\HH^{\otimes m}$,}$$
where
 $$P^{(n)}:=\sum_{\sigma\in\SSn}\ff(\sigma)\qquad\qquad
\text{(note $P^{(0)}:=\id$)}$$
is the canonic operator corresponding to the quasi-multiplicative
function $\ff:\SSn\to B(\HH^{\otimes n})$ given by $\ff(e)=\id$ and
$\ff(\pi_i)=T_i$ ($i=1,\dots,n-1$). According to Theorem 2.2 the
operators $P^{(n)}$ are positive, thus $\langle\enskip,\enskip\rangle_T$
is positive definite. If $\Vert T\Vert<1$ then, by 2.3, we know that all
$P^{(n)}$ are strictly positive and we can take as $\cFT$ the completion
of $\cF$ with respect to $\langle\enskip,\enskip\rangle_T$. In the case
$\Vert T\Vert=1$, we might get a kernel of
$\langle\enskip,\enskip\rangle_T$ and we have to divide this out before
taking the completion.
\proclaim{Theorem 3.1} i) For all $f\in\HH$, $d(f)$ and $d^*(f)$ are
adjoints of each other on $\cFT$, i.e. for all $k\in\NN$ and all
$\xi,\eta\in\bigoplus_{n=0}^k\HH^{\otimes n}$ we have
 $$\la d^*(f)\xi,\eta\raT=\la\xi,d(f)\eta\raT.$$\noindent
ii) If $\Vert T\Vert=q<1$, then
 $$\Vert d^*(f)\Vert_T\leq \Vert f\Vert \frac 1{\sqrt{1-q}}.$$
\endproclaim
\demo{Proof}
i)
By definition of the $T_i$, we have
 $$l^*(f)T_i=T_{i+1}l^*(f)\qquad (i\geq 1),$$
which implies
 $$l^*(f)P^{(n)}=(\id\otimes P^{(n)})l^*(f) \quad\text{or}\quad
\Pn l(f)=l(f)(\id\otimes\Pn).$$
Furthermore, our general decomposition
 $P(W)=P(D_J)P(W_J)$
gives for the case
 $$W=\SS_{n+1}, \quad J=\{\pi_2,\pi_3,\dots,\pi_n\},\quad
D_J=\{e,\pi_1,\pi_2\pi_1,\dots,\pi_{n-1}\pi_{n-2}\dots\pi_1\}$$
the relation
 $$P^{(n+1)}=P(\SS_{n+1})=R^{(n+1)*}(\id\otimes P^{(n)})=
(\id\otimes\Pn)R^{(n+1)},$$
where
 $$R^{(n)}:=\id+T_1+T_1T_2+\dots+T_1\dots T_{n-2}T_{n-1}.$$
Note that
 $$d(f)=l(f)R^{(n)}\qquad\text{on $\HH^{\otimes n}$.}$$
We have now for $\xi\in\HH^{\otimes n}$ and $\eta\in\HH^{\otimes (n+1)}$
 $$\align
\la d^*(f)\xi,\eta\raT&=\la d^*(f)\xi,P^{(n+1)}\eta\ra\\
&=\la \xi,l(f) P^{(n+1)}\eta\ra\\
&=\la\xi,l(f)(\id\otimes\Pn)R^{(n+1)}\eta\ra\\
&=\la\xi,\Pn l(f) R^{(n+1)}\eta\ra\\
&=\la\xi,P^{(n)}d(f)\eta\ra\\
&=\la \xi,d(f)\eta\raT.\endalign$$
ii) Since
 $$\Vert R^{(n)}\Vert\leq 1+q+q^2+\dots+ q^{n-1}\leq \frac 1{1-q},$$
we have
 $$\align
P^{(n+1)}P^{(n+1)}&=(\id\otimes P^{(n)})R^{(n+1)}R^{(n+1)*}(\id\otimes\Pn)\\
&\leq \frac 1{(1-q)^2}(\id\otimes \Pn)(\id\otimes\Pn)\endalign$$
hence (because of the positivity of $P^{(n)}$ and $P^{(n+1)}$)
 $$P^{(n+1)}\leq \frac 1{1-q} (\id\otimes\Pn),$$
which gives for $\xi\in\HH^{\otimes n}$
 $$\align
\Vert d^*(f)\xi\Vert_T^2&=\la d^*(f)\xi,d^*(f)\xi\raT\\
&=\la f\otimes\xi,f\otimes\xi\raT\\
&=\la f\otimes\xi,P^{(n+1)}f\otimes\xi\ra\\
&\leq\frac 1{1-q}\la f\otimes\xi,(\id\otimes\Pn)f\otimes\xi\ra\\
&=\frac 1{1-q}\la f,f\ra\la \xi,\Pn\xi\ra\\
&=\frac 1{1-q}\Vert f\Vert^2 \thinspace\Vert \xi\Vert_T^2.\endalign$$
\line{\hfill $\diamondsuit$}
\enddemo
If we choose some basis $\{e_i\}_{i\in I}$ of $\HH$ and define the matrix
$t$ by
 $$T e_a\otimes e_b=\sum_{c,d\in I} t_{ab}^{dc} e_d\otimes e_c \qquad
(a,b\in I),$$
then, by using the definition of our creation and annihilation
operators, it is easy to check that the operators $d_i:=d(e_i)$ ($i\in I$)
fulfill the relations
 $$d_id_j^*-\sum_{r,s\in I} t_{js}^{ir} d_r^* d_s=\delta_{ij}\qquad
(i,j\in I).$$
Since by construction $d(f)\Omega=0$ for all $f\in\HH$, we have obtained
the Fock representation of these relations.\par
Now we want to recover the $q_{ij}$-relations from our general construction.
For this we consider the
operator $T=Q\pi_1$, where $Q$ is the multiplication operator
 $$Q(e_i\otimes e_j)=q_{ij}(e_i\otimes e_j)$$
and $\pi_1$ the natural action of the corresponding transposition
 $$\pi_1(f\otimes g)=g\otimes f \qquad (f,g\in\HH).$$
This $T$ is self-adjoint (because of $\bar q_{ij}=q_{ji}$), contractive
($\Vert T\Vert=\sup_{i,j\in I}\vert q_{ij}\vert\leq 1$), and fulfills the
braid relation (BR).
Thus the foregoing construction may be applied to it. In
this case one gets the following concrete formula for the annihilation
operator.
 $$\multline
d(e_i)e_{i(1)}\otimes\dots\otimes e_{i(n)}=\\
=\sum_{k=1}^n q_{i(k),i(k-1)}q_{i(k),i(k-2)}\dots q_{i(k),i(1)}
\delta_{i,i(k)} e_{i(1)}\otimes\dots\otimes\check e_{i(k)}\otimes\dots
\otimes e_{i(n)},\endmultline
 $$
where $e_{i(k)}$ has to be deleted in the tensor.\par
Thus we get the following corollary on the existence of Fock
representations out of our constructions.
\proclaim{Corollary 3.2}
i) Let $T\in B(\HH\otimes\HH)$ be a self-adjoint contraction
fulfilling the braid relation and write
 $$Te_a\otimes e_b=\sum_{c,d\in I} t_{ab}^{dc} e_d\otimes e_c\qquad
(a,b\in I)$$
for some basis $\{e_i\}_{i\in I}$ of $\HH$. Then there exist operators
$d_i$ ($i\in I$) on some Hilbert space $\tilde\HH$ and a \lq vacuum
vector' $\Omega\in\tilde\HH$ such that $d_i\Omega=0$ for all $i\in I$
and
 $$d_id_j^*-\sum_{r,s\in I} t_{js}^{ir} d_r^* d_s=\delta_{ij}\qquad
(i,j\in I).$$
If $\Vert T\Vert<1$, then the $d_i$ are bounded.\newline
ii) In particular,
for given $q_{ij}$ ($i,j\in I$) with $\bar q_{ij}=q_{ji}$ for all $i,j\in I$
and $\sup_{i,j\in I} \vert q_{ij}\vert=q\leq 1$
there exist operators $d_i$ ($i\in I$) on some Hilbert space $\tilde\HH$
and a \lq vacuum vector' $\Omega\in\tilde\HH$ such that $d_i\Omega=0$ for
all $i\in I$ and
 $$d_id_j^*-q_{ij}d_j^*d_i=\delta_{ij}\id \qquad (i,j\in I).$$
If $q<1$, then the $d_i$ are bounded.
\endproclaim
\demo{Remarks}
1) Consider the $q_{ij}$-relations.
Let $q<1$. Then, for $q_{ii}\geq0$, it follows from
 $d_id_i^*-\qii d_i^*d_i=\id$ the estimate $\Vert d_i\Vert^2\leq 1+\qii
\Vert d_i\Vert^2$, i.e. $\Vert d_i\Vert\leq 1/\sqrt{1-\qii}$. For
$\qii<0$, we even have $d_id_i^*=\id+\qii d_i^*d_i\leq\id$, i.e.
$\Vert d_i\Vert\leq1$. The restriction of our representation from
$\cF_T$ to the linear span of $\{e_i^{\otimes n}\mid n\in\NN\}$
shows that these inequalities are indeed equalities, thus
 $$\Vert d_i\Vert_T=\cases \frac 1{\sqrt{1-\qii}},&\text{if $\qii\geq0$}\\
1,&\text{if $\qii<0$}.\endcases$$
This is true for all $q<1$, thus, by continuity, also for $q=1$.\newline
2) For the crucial step, namely the positivity of all $\Pn$,
in our construction of the $q_{ij}$-relations
we do not need 2.2 in full generality but only for the special
case of $T_i=Q_i\pi_i$, where the $Q_i$ commute. For this case a simpler
proof of 2.2 (for $W=\SSn$) is given in \l JSW2\r\  (but without any assertion
on strict positivity of $P^{(n)}$ in the case $\Vert T\Vert<1$).
\enddemo
\vskip1cm
\heading
{\bf 4. Operator spaces}
\endheading
Now we want to consider the deformed commutation relations constructed in
the last section from an operator space point of view. Operator spaces
were introduced by Effros and Ruan \l ER1,ER2\r\ and further investigated
by Blecher and Paulsen \l BP\r\ and
Pisier \l Pis2\r. Operator spaces are closed linear subsets of $B(\HH)$ for
some Hilbert space $\HH$ and have a lot of nice properties. The
philosophy behind their introduction is that they quantize functional
analysis in that sense that in the usual statements, e.g. in norm
inequalities, numbers are replaced by operators. We refer to \l ER2,BP,Pis2\r\
for more details.\par
One canonic operator space is the Hilbert space $R\cap C\subset M_N\oplus
M_N$ (where $M_N$ are the $N\times N$-matrices, for $N=\infty$ the compact
operators)
with basis $\{\delta_i\}_{i=1,\dots,N}$ given by
 $$\delta_i=
\pmatrix 0&\qquad&\qquad&\qquad\\
\vdots&\qquad&\qquad&\qquad\\
1&\qquad&\bigcirc&\qquad\\
\vdots&\qquad&\qquad&\qquad\\
0&\qquad&\qquad&\qquad\endpmatrix\oplus
\pmatrix
0&\hdots&1&\hdots&0\\
\quad&\quad&\quad&\quad\quad\\
\quad&\quad&\quad&\quad\quad\\
\quad&\quad&\bigcirc&\quad&\quad\\
\quad&\quad&\quad&\quad\quad\\
\quad&\quad&\quad&\quad\quad
\endpmatrix,$$
where the 1 is appearing in the $i$-th position in the first column or
first row, respectively. Operator spaces which are also Hilbert spaces
are called Hilbertian operator spaces.
The Hilbertian operator space $R\cap C$
has the following characterizing property: For all $a_i\in B(\tilde\HH)$
($i=1,\dots,N$) on some Hilbert space $\tilde\HH$ one has
 $$\Vert \sum_{i=1}^N a_i\otimes \delta_i\Vert_{B(\tilde\HH)\otimes M_N}
=\Vert (a_1,\dots,a_N)\Vert_{\max},$$
where
 $$\Vert(a_1,\dots,a_N)\Vert_{\max}:=\max\bigl(
\Vert\sum_{i=1}^N a_ia_i^*\Vert^{1/2},\Vert\sum_{i=1}^N a_i^*a_i\Vert^{1/2}
\bigr).$$
We consider now the operators $d(f)$ and $d^*(f)$ ($f\in\HH$) on $\cF_T$
as constructed, for a given self-adjoint contraction $T$ fulfilling the
braid relation (BR),
in the last section. Assume in the following $\Vert T\Vert=
q<1$. We choose a basis $\{e_i\}$ of $\HH$ and put $d_i:=d(e_i)$ and
 $$G_i:=d_i+d_i^*.$$
Then we claim that the operator space generated by the closure of
the linear span of the $G_i$ is, as an operator space, isomorphic
to $R\cap C$,
where $N=\dim \HH$. This means nothing else than the following
norm estimate.
\proclaim{Theorem 4.1}
We have for arbitrary operators $a_i\in B(\tilde\HH)$ ($i=1,\dots,N$)
with $N\leq \dim\HH$ the estimate
 $$\Vert(a_1,\dots,a_N)\Vert_{\max}\leq
\Vert\sum_{i=1}^N a_i\otimes G_i\Vert_{\HsFT}\leq\frac 2{\sqrt{1-q}}
\Vert(a_1,\dots,a_N)\Vert_{\max}.$$
\endproclaim
The case $T=0$ was treated by Haagerup and Pisier \l HP\r. Our proof will
follow their ideas.
\demo{Proof}
By definition $d(f)=l(f)R^{(n)}$ on $\HH^{\otimes n}$, i.e. more
generally $d(f)=l(f)R$, where the infinite sum
 $$R:=\id+T_1+T_1T_2+T_1T_2T_3+\dots$$
makes sense because of $\Vert T\Vert <1$. The crucial step is now a norm
estimate for this $R$. Of course, we have with respect to the usual norm
on the full Fock space $\Vert R\Vert\leq 1/(1-q)$. We want to show
that the same estimate is true for $\Vert R\Vert_T$. Remember that
$\Pnpe=(\id\otimes\Pn)\Rnpe$ and
 $\Pnpe\leq 1/(1-q)(1\otimes\Pn)$.
Then we have for $\xi\in\HH^{\otimes(n+1)}$
 $$\align
\Vert R\xi\Vert^2_T&=\la \Rnpe\xi,\Rnpe\xi\raT\\
&=\la\Rnpe\xi,\Pnpe\Rnpe\xi\ra\\
&\leq \frac 1{1-q}\la\Rnpe\xi,(\id\otimes\Pn)\Rnpe\xi\ra\\
&=\frac 1{1-q}\la\Rnpe\xi,\Pnpe\xi\ra\\
&=\frac 1{1-q}\la\Rnpe\xi,\xi\raT\\
&\leq\frac 1{1-q}\Vert R\xi\Vert_T\thinspace \Vert\xi\Vert_T,\endalign$$
which implies
 $$\Vert R\xi\Vert_T\leq\frac 1{1-q}\Vert\xi\Vert_T,\qquad
\text{hence}\qquad \Vert R\Vert_T\leq \frac 1{1-q}.$$
Since
 $$\Vert \sum_{i=1}^N b_ic_i\Vert\leq\Vert\sum_{i=1}^N b_ib_i^*\Vert^{1/2}
\thinspace \Vert\sum_{i=1}^Nc_i^*c_i\Vert^{1/2}$$
for arbitrary bounded operators $b_1,\dots,b_N,c_1,\dots,c_N$ on some
Hilbert space, we obtain now
 $$\align
\Vert\sum_{i=1}^N a_i\otimes d_i^*\Vert_{\HsFT}&=
\Vert\sum_{i=1}^N (\id\otimes d_i^*)(a_i\otimes\id)\Vert_{\HsFT}\\
&\leq\Vert\sum_{i=1}^N d_i^*d_i\Vert^{1/2}_T\thinspace
\Vert\sum_{i=1}^N a_i^*a_i\Vert_{\Hs}^{1/2}.\endalign$$
Because of
 $$\sum_{i=1}^N d_i^*d_i=(\sum_{i=1}^N l_i^*l_i)R=(1-P_\Omega)R,$$
where $P_\Omega$ is the projection onto the vacuum $\Omega$, we get
 $$\align
\Vert\sum_{i=1}^N a_i\otimes d_i^*\Vert_{\HsFT}&\leq \Vert R\Vert_T^{1/2}
\thinspace\Vert\sum_{i=1}^N a_i^*a_i\Vert_{\Hs}^{1/2}\\
&\leq\frac 1{\sqrt{1-q}}\Vert\sum_{i=1}^N a_i^*a_i\Vert_{\Hs}^{1/2}\endalign$$
and by taking adjoints
 $$\Vert\sum_{i=1}^N a_i\otimes d_i\Vert_{\HsFT}\leq
\frac 1{\sqrt{1-q}}\Vert\sum_{i=1}^N a_ia_i^*\Vert_{\Hs}^{1/2},$$
which yield together the right inequality of our assertion.\par
For the other inequality we use the vacuum expectation state
 $$\ee(T):=\la\Omega,T\Omega\ra \qquad (T\in B(\cF_T)).$$
We only need
 $$\ee(G_iG_j)=\la\Omega,d_id_j^*\Omega\ra=\delta_{ij}$$
to obtain
 $$\align
\Vert\sum_{i=1}^N a_i\otimes G_i\Vert^2_{\HsFT}&\geq
\sup\Sb \text{$\ff$ state}\\ \text{on $B(\Hs)$}\endSb (\ff\otimes\ee)
\bigl\lb (\sum_{i=1}^N a_i\otimes G_i)^*(\sum_{j=1}^N a_j\otimes G_j)\bigr\rb\\
&=
\sup\Sb \text{$\ff$ state}\\ \text{on $B(\Hs)$}\endSb \ff(\sum_{i=1}^N
a_i^*a_i)\\
&=\Vert\sum_{i=1}^N a_i^*a_i\Vert,\endalign$$
and analogously
 $$\Vert\sum_{i=1}^N a_i\otimes
G_i\Vert^2_{\HsFT}\geq \Vert \sum_{i=1}^N a_ia_i^*\Vert.$$
\line{\hfill $\diamondsuit$}
\enddemo
Our theorem characterizes completely the operator space structure of our
deformations, namely this structure does not depend on the deformation
(at least as long as $\Vert T\Vert<1$). One may also ask about the
$C^*$- or $W^*$-structure of our deformations. In this respect, the
situation is not so clear. Let us make in the following
some remarks in this direction.\par
For $\vert q\vert$ sufficiently small, the method of \l JSW1\r\  should
still work showing that the $C^*$-algebra generated by all $d(f)$ ($f\in\HH$)
is equal to the extension of the Cuntz algebra $O_n$ by the compact operators,
where $n=\dim\HH$. See \l JSW2\r\ for investigations in this direction. It is
conceivable that the $C^*$-structure of the $q_{ij}$-relations or even
of our general deformations is the same for
all $\Vert T\Vert<1$.\par
Another interesting problem is the structure of the von Neumann algebra
$\cM_T$ generated by all $G_i$. Typically, these von Neumann algebras are
not injective. Injectivity of a von Neumann algebra $\cM\subseteq B(\HH)$
means that there exists a projection of norm 1 from $B(\HH)$ onto $\cM$.
\proclaim{Theorem 4.2}
If the vacuum expectation $\ee$ is a faithful trace on $\cM_T$ and
$\dim\HH>16/(1-q)^2$, then $\cM_T$ is not injective.
\endproclaim
\demo{Proof}
If $\cM_T$ were injective we would have for all $a_i,b_i\in\cM_T$
(compare \newline Corollary 2 of \l Was\r)
 $$\Vert\sum_{i=1}^m a_i\otimes
\bar b_i\Vert \geq\vert \ee(\sum_{i=1}^m a_i b_i^*)\vert
, $$
in particular, for $a_i=b_i=G_i$,
 $$\Vert\sum_{i=1}^m G_i\otimes
\bar G_i\Vert\geq \ee(\sum_{i=1}^m G_iG_i)=m.$$
But on the other side, by putting
$a_i=\bar G_i$ in Theorem 4.1, we also have
 $$\Vert\sum_{i=1}^m G_i\otimes\bar G_i\Vert\leq \frac 4{1-q}\sqrt m,$$
which leads, for $m>16/(1-q)^2$, to a contradiction.\newline
\line{\hfill $\diamondsuit$}
\enddemo
By following the ideas of Theorem 2.9 in \l Pis2\r, we see that
faithfulness of $\ee$ is not really needed. But according to the next
theorem we do not need to make this distinction.
\proclaim{Theorem 4.3}
Assume that the vacuum expectation $\ee$ is tracial on $\cM_T$. Then the
vacuum $\Omega$ is cyclic and separating for $\cM_T$. In particular, $\ee$
is faithful.
\endproclaim
\demo{Proof}
For cyclicity of $\Omega$ it suffices to see that we can obtain all
basis product vectors $e_{i(1)}\otimes\dots\otimes e_{i(k)}$ for all
$k\in\NN$ and all $i(1),\dots,i(k)\in\{1,\dots,n\}$ from $\Omega$ by
application of some polynomial in the $G_i$. Since
 $$e_{i(1)}=G_{i(1)}\Omega$$
and
 $$e_{i(1)}\otimes\dots\otimes e_{i(k)}=G_{i(1)}\dots G_{i(k)}\Omega
-\eta \qquad\text{with}\quad \eta\in\bigoplus_{l=0}^{k-1}\HH^{\otimes l},$$
this follows by induction.\newline
To show that $\Omega$ is separating for $\cM_T$ is the same as showing that
$\Omega$ is cyclic for $\cM_T'$. Let us define the anti-linear conjugation
operator
$J:\cF_T\to\cF_T$ by $JA\Omega=A^*\Omega$ for $A\in\cM_T$. This is
well-defined because the trace property of $\ee$ implies
 $\Vert A\Omega\Vert_{\cF_T}=\Vert A^*\Omega\Vert_{\cF_T}$.
One can easily check that
 $J\cM_T J\subseteq \cM_T'$.
Since $\Omega$ is cyclic for $J\cM_T J$, the assertion follows.\newline
\line{\hfill $\diamondsuit$}
Note that we have shown that $\cM_T$ is in standard form and thus
$J\cM_T J=\cM_T'$.
By the way, the conjugation operator $J$ is explicitly given by
 $$J(e_{i(1)}\otimes\dots\otimes e_{i(n)})=e_{i(n)}\otimes\dots\otimes
e_{i(1)},$$
i.e.
the operators $JG_iJ$ are like the $G_i$, only
action from the left is replaced by action from the right.\par
This raises the question whether $\ee$ is a trace on $\cM_T$.
This can be answered like follows.
\proclaim{Theorem 4.4}
1) The vacuum expectation $\ee$ is a trace on $\cM_T$ if and only if $T$
fulfills
 $$\la e_d\otimes e_c,T e_a\otimes e_b\ra=
\la e_c\otimes e_b,T e_d\otimes e_a\ra$$
for all $a,b,c,d\in\{1,\dots,\dim\HH\}$.\newline
2) In particular, in
the case of the $q_{ij}$-relations,
$\ee$ is a trace if and only if the $q_{ij}$ are
symmetric and hence real.
\endproclaim
\demo{Proof}
1) We will only give a sketch of the proof for the general case.\newline
Let us write
 $$T e_a\otimes e_b=\sum_{c,d} t_{ab}^{dc} e_d\otimes e_c.$$
Then the asserted condition for $T$ reads as
 $$t_{ab}^{dc}=t_{da}^{cb},$$
i.e. $t_{ab}^{dc}$ is cyclic in its four indices. (By the way, $T=T^*$ reads
in this language as $\bar t_{ab}^{dc}=t_{dc}^{ab}$.)\newline
Necessity of this cyclicity condition follows from
 $$\ee(G_aG_bG_cG_d)=\ee(d_ad_bd^*_cd^*_d)+\ee(d_ad_b^*d_cd_d^*)=
\la e_a,e_d\ra\thinspace \la e_b,e_c\ra+ t_{cd}^{ba}+
\la e_a,e_b\ra\thinspace \la e_c,e_d\ra.$$
To see that cyclicity is also a sufficient condition, one has to check,
by using the very definition of $d_i$ and $d_i^*$, the following
formula.
 $$\ee(G_{i(1)}\dots G_{i(m)})=\cases 0,&\text{$m$ odd}\\
\sum_{\cV}\cT(\cV),&m=2r,\endcases$$
where the sum runs over all pairings $\cV=\{(a_1,z_1),\dots,(a_r,z_r)\}$
of the indices $1,\dots,2r$ (we always assume $a_k<z_k$ and $a_k<a_l$ for
$k<l$) and where $\cT(\cV)$ is a factor which is calculated from a given
$\cV$ in the following way: Put $2r$ points on a circle and denote them
in clock-wise order by $1,\dots,2r$. Connect the points $a_k$ and $z_k$
for all $k=1,\dots,r$ by an arc inside the circle in such a way that at most
two arcs cross in one point and such that the number of these crossing
points is minimal. Thus we get a graph consisting of points, namely
the outer points on the circle and the crossing points, and edges,
namely the pieces of our arcs connecting two points. To each edge, we
assign a variable $a,b,c,\dots$. This graph determines now $\cT(\cV)$
by the following rules: Each outer point $k$ with edge $a$ gives a
factor $\delta_{i(k),a}$. Each crossing point with the four edges
$a,b,c,d$ (in clock-wise order) gives a factor $t_{ab}^{dc}$. Take
then the product over the factors corresponding to all points and sum this
over all variables of the edges, each running from 1 to $\dim\HH$. The
result is $\cT(\cV)$.\newline
Examples: For $\cV=\{(1,4),(2,4)\}$
we have
 $$\cT(\cV)=\sum_{a,b}\delta_{a,i(1)}\delta_{a,i(4)}\delta_{b,i(2)}\delta_
{b,i(3)}=\delta_{i(1),i(4)}\delta_{i(2),i(3)}.$$
For $\cV=\{(1,3),(2,5),(4,6)\}$ we have
 $$\cT(\cV)=\sum_a t_{i(2)i(3)}^{i(1)a} t_{i(4)i(5)}^{ai(6)},$$
whereas for $\cV=\{(1,4),(2,5),(3,6)\}$ we obtain
 $$\cT(\cV)=\sum_{a,b,c} t_{i(6)i(1)}^{cb} t_{i(2)i(3)}^{ba} t_{i(4)i(5)}^
{ac}.$$
Note that the braid relation for $T$ ensures that $\cT(\cV)$ does not
depend on the way we have drawn our graph, as long as we keep the
number of crossing points minimal.
If we do not assume cyclicity of $t_{ab}^{dc}$,
then a similar formula would be valid, one only has to take care how to
arrange the variables at $t_{ab}^{dc}$ at the crossing points. For this
one has to distinguish between ingoing and outgoing edges.\newline
Having the above formula for the calculation of $\ee$, one sees quite
easily that under a cyclic permutation of $\cV$ the clockwise order
at the crossing points does not change, thus $\cT(\cV)$ does not
change under such a cyclic permutation (under the assumption
$t_{ab}^{dc}=t_{da}^{cb}$) and hence $\ee$ is a trace on $\cM_T$.\newline
2) Since
 $t_{ab}^{dc}=q_{ba}\delta_{bd}\delta_{ca}$, we have
$t_{ab}^{dc}=t_{da}^{cb}$ if and only if $q_{ab}=q_{ba}$.\newline
In this case, $\cT(\cV)$ from part 1 can be written more explicitly
as
 $$\cT(\cV)=
\delta_{i(a_1),i(z_1)}\dots\delta_{i(a_r),i(z_r)}\cdot
\ct(\cV)$$
for a given pairing $\cV=\{(a_1,z_1),\dots,(a_r,z_r)\}$
of the indices $1,\dots,2r$. The number
 $\ct(\cV)$ denotes
a weighting factor taking into account the number of inversions
of $\cV$, namely with
 $$I(\cV):=\{(k,l)\mid k,l=1,\dots,r\quad\text{such that}\quad
a_k<a_l<z_k<z_l\}$$
it is given by
 $$\ct(\cV)=\prod_{(k,l)\in I(\cV)} q_{i(a_k),i(a_l)}.$$
In this case the
formula for $\ee(G_{i(1)}\dots G_{i(m)})$ can be proved quite easily from
the identical one for
$\ee(d_{i(1)}^{\#}\dots d_{i(m)}^{\#})$, where
$d_{i(k)}^{\#}$ stands for $d_{i(k)}$ or $d_{i(k)}^*$. This latter formula
follows by noticing that it is true for products of the form \newline
$d_{i(1)}^*\dots d_{i(k)}^*d_{i(k+1)}\dots d_{i(k+l)}$ and that both
sides of the formula change in the same way under application of the
$q_{ij}$-relations, see \l BSp1,Spe2\r.
\newline
\line{\hfill $\diamondsuit$}
\enddemo
We conjecture that $\cM_T$ is, at least for the $q_{ij}$-relations,
a factor. This will be pursued further in forthcoming investigations.
\enddemo
\vskip1cm
\heading
{\bf 5. Completely positive maps corresponding to block length}
\endheading
The completely positive maps on Coxeter groups considered in Sect. 2 were
canonical generalizations of the basic example
 $\ff(\sigma)=q^{\vert\sigma\vert}$,
where $\vert\sigma\vert$ is the usual length function on our Coxeter
group $W$. This example appeared (for $W=\SSn$) quite naturally in the
course of our investigations on generalized Brownian motions in
\l BSp1\r. In \l BSp3\r\ we considered another example of a Brownian
motion which is intimately connected with Voiculescu's concept of
freeness \l VDN\r.
We found that once more the positive definiteness of some function
on $\SSn$ is the key point in our construction. This function is again of
the form
 $$\ff(\sigma)=q^{\Vert \sigma\Vert},$$
but now $\Vert\sigma\Vert$ is another length function on $\SSn$. Namely,
whereas $\vert\sigma\vert$ counts the number of generators in a reduced
representation of $\sigma$, the function $\Vert\sigma\Vert$ gives the
number of {\it different} generators.
This length function and the corresponding quasi-multiplicative
$\ff$ can now again be extended in a canonical way to arbitrary Coxeter
groups and to operator valued functions.\par
Let $(W,S)$ be an arbitrary Coxeter group. If $\sigma=s_{i(1)}\dots
s_{i(k)}$ is a reduced representation of $\sigma$, then we put
 $$b(\sigma):=\{s_{i(1)},\dots,s_{i(k)}\},$$
the set of generators appearing in $\sigma$. Although a reduced
representation is not unique, $b(\sigma)$ is well defined, see \l Bou\r.
For example, in $W=\SS_3$, we have $\pi_1\pi_2\pi_1=\pi_2\pi_1\pi_2$,
and $b(\pi_1\pi_2\pi_1)=b(\pi_2\pi_1\pi_2)=\{\pi_1,\pi_2\}$.\par
We call the corresponding length function
 $$\Vert\sigma\Vert:=\#b(\sigma)$$
block length function. As will follow from our Theorem 5.1, it is a
positive definite function on $W$. For $W=\SS_n$, it has a nice
graphical meaning, namely $$\Vert\sigma\Vert=
n-\text{the number of connected components of the graph of $\sigma$},$$
e.g.
 $$\sigma=\pi_1\pi_4\pi_3=\matrix \underset \circ\to 1&
                  \underset \circ\to 2&
                  \underset \circ\to 3&
                  \underset \circ\to 4&
                  \underset \circ\to 5\\
\thinspace\\
\overset\circ\to 1&
 \overset \circ\to 2&
 \overset \circ\to 3&
 \overset \circ\to 4&
 \overset \circ\to 5
\endmatrix\qquad\qquad
\Vert\sigma\Vert=5-2=3$$
 $$\sigma=e=
\matrix \underset \circ\to 1&
                  \underset \circ\to 2&
                  \underset \circ\to 3&
                  \underset \circ\to 4&
                  \underset \circ\to 5\\
\thinspace\\
 \overset \circ\to 1&
 \overset \circ\to 2&
 \overset \circ\to 3&
 \overset \circ\to 4&
 \overset \circ\to 5
\endmatrix\qquad\quad\Vert\sigma\Vert=5-5=0.$$
The analogue of 2.1 for this concept of block length is now the
following.
\proclaim{Theorem 5.1}
Let $T_i\in B(\HH)$ ($i=1,\dots,n$) be bounded operators on some Hilbert
space $\HH$, which satisfy:\roster
\item"i)"
 $0\leq T_i\leq\id$ for all $i=1,\dots,n$.
\item"ii)"
The $T_i$ commute:
 $T_iT_j=T_jT_i$ for all $i,j=1,\dots,n$.
\endroster
Define now a quasi-multiplicative (with respect to $\Vert\sigma\Vert$)
function
 $$\ff:\CC W\to B(\HH) \qquad\text{by}\qquad
\ff(\sigma):=\prod\Sb \text{$i$ with}\\
s_i\in b(\sigma)\endSb T_i\qquad (\ff(e):=\id).$$
Then $\ff$ is completely positive.
\endproclaim
\demo{Remarks}
Note that our assumptions on the $T_i$ are quite natural.\newline
1) In the example
$\ff_q(\sigma)=q^{\Vert\sigma\Vert}$ in the case $W=\SSn$ one can
check that $\ff_q$ is positive definite only for $1\geq q\geq\alpha_n$,
where $\alpha_n<0$, but $\lim_{n\to\infty}\alpha_n=0$. Thus, in general,
we have to assume $T_i\geq 0$.\newline
2) Also commutativity of the $T_i$ is necessary, otherwise the relations
in $W$ would conflict with a canonical definition of $\ff$, e.g. for
$W=\SS_3$ and $\pi:=\pi_1\pi_2\pi_1=\pi_2\pi_1\pi_2$ there is no
canonic preference for one of the two possibilities $\ff(\pi)=T_1T_2$ or
$\ff(\pi)=T_2T_1$, thus they should coincide.
\enddemo
\demo{Proof}
Since the pointwise product of two commuting completely positive maps is
again completely positive \l Boz3\r, it suffices to consider the special
case where all $T_i$ but one are equal to $\id$, i.e. for arbitrary but
fixed $k\in\{1,\dots,n\}$ we consider
 $$\ff(\sigma)=\cases T,&\text{if $s_k\in b(\sigma)$}\\
                      \id,&\text{if $s_k\not\in b(\sigma)$,}\endcases$$
where $T:=T_k$ fulfills $0\leq T\leq\id.$\newline
Since $T$ can be diagonalized by the spectral theorem, the assertion can
be reduced to the scalar valued case and we only have to treat the special
cases
 $$\ff_q(\sigma)=\cases q,&\text{if $s_k\in b(\sigma)$}\\
                      1,&\text{if $s_k\not \in b(\sigma)$}\endcases$$
for all $k\in\{1,\dots,n\}$ and all $q$ with $0\leq q\leq 1$.\newline
Let $W_k$ be the parabolic subgroup of $(W,S)$ generated by
$J:=S\backslash\{s_k\}$. Then one knows \l Bou\r\ that
 $$s_k\in b(\sigma)\Longleftrightarrow \sigma\cdot W_k\not= W_k.$$
Now consider the kernel $\delta$ on all subsets of $W$ given by
 $$\delta(A,B)=\cases 1,&\text{if $A=B$}\\
0,&\text{if $A\not= B$}\endcases\qquad\text{for $A,B\subseteq W$.}$$
Then we have by putting $q=\exp(-t)$ ($0<t<\infty$)
 $$\ff_q(\sigma)=q^{1-\delta(\sigma\cdot W_k,W_k)}=
e^{-t}e^{t\delta(\sigma\cdot W_k,W_k)}$$
or
 $$\ff_q(\tau^{-1}\sigma)=e^{-t}e^{t\delta(\sigma\cdot W_k,\tau\cdot W_k)}.$$
Since $\delta$ is positive definite on all subsets of $W$ we get, by the
Sch\"onberg theorem (see, e.g., \l Boz3\r), the positive definiteness of
$\ff_q$ for $t>0$. The case $q=1$ is trivial, and $q=0$ follows by
continuity from $q\searrow 0$.\newline
\line{\hfill $\diamondsuit$}
\enddemo
\demo{Remarks}
1) Note that, contrary to the situation considered in Sect. 2, the
scalar valued case contains all essential information, the operator
valued version is a mere transcription to diagonal operators. Thus, in
the spirit of the remarks 2 and 3 at the end of Sect. 2, we are not
restricted to amenable Coxeter groups, but Theorem 5.1 is valid for all
Coxeter groups.\newline
2) Note that here there is no reduction to a positivity problem for some
operator $P$ like the reduction from 2.1 to 2.2. The statement that
$\sum_{\sigma\in W}\ff(\sigma)\geq 0$,
is trivially true because of $T_i\geq 0$ and $T_iT_j=T_jT_i$, but it is
by far not sufficient for the complete positivity of $\ff$.
\enddemo
\vskip1cm
\heading
{\bf Acknowledgements}
\endheading
This work has been supported by grants from KBN (M.B.) and the Deutsche
Forschungsgemeinschaft (R.S.).\par
The first author thanks G. Pisier for a nice working atmosphere during
his stay at Paris in March - April 93.
\vskip1cm
\heading
{\bf References}
\endheading
\roster
\item"\l BP\r"
Blecher, D., Paulsen, V.: Tensor products of operator spaces.
J. Funct. Anal. {\bf 99}, 262-292 (1991)
\item"\l Boz1\r"
Bo\D zejko, M.: Positive definite kernels, length functions on groups and a
non commutative von Neumann inequality. Studia Math. {\bf 95},
107-118 (1989)
\item"\l Boz2\r"
Bo\D zejko, M.: In peparation
\item"\l Boz3\r"
Bo\D zejko, M.: Positive and negative definite kernels on discrete groups.
To appear in {\it Springer Lecture Notes}
\item"\l BJS\r"
Bo\D zejko, M., Januszkiewicz, T., Spatzier, R.J.: Infinite Coxeter
groups do not have Kazhdan's property. J. Operator Theory
{\bf 19}, 63-68 (1988)
\item"\l BSp1\r"
Bo\D zejko, M., Speicher, R.: An Example of a Generalized Brownian
Motion. Comm. Math. Phys. {\bf 137}, 519-531 (1991)
\item"\l BSp2\r"
Bo\D zejko, M., Speicher, R.: An example of a generalized Brownian
motion II. In: Quantum Probability and Related Topics VII
(ed. L. Accardi, pp. 67-77), Singapore: World Scientific 1992
\item"\l BSp3\r"
Bo\D zejko, M., Speicher, R.: Interpolations between bosonic and fermionic
relations given by generalized Brownian motions. Preprint SFB 123-691,
Heidelberg 1992
\item"\l Bou\r"
Bourbaki, N.: Groupes et algebres de Lie, Chap. 4,5,6. Paris: Hermann 1968
\item"\l Car\r"
Carter, R.W.: Simple groups of Lie type. New York: John Wiley \& Sons 1972
\item"\l deH\r"
de la Harpe, P.: Groupes de Coxeter non affines. Exposition Math.
{\bf 5}, 91-95 (1987)
\item"\l ER1\r"
Effros, E., Ruan, Z.J.: On matricially normed spaces. Pacific J.
Math. {\bf 132}, 243-264 (1988)
\item"\l ER2\r"
Effros, E., Ruan, Z.J.: A new approach to operator spaces. Canadian
Math. Bull. {\bf 34}, 329-337 (1991)
\item"\l Eva\r"
Evans, D.E.: On $O_n$. Publ. RIMS {\bf 16}, 915-927 (1980)
\item"\l Fiv\r"
Fivel, D.: Interpolation between Fermi and Bose statistics using
generalized commutators. Phys. Rev. Lett. {\bf 65}, 3361-3364 (1990)
\item"\l Gre\r"
Greenberg, O.W.: Particles with small violations of Fermi or Bose statistics.
Phys. Rev. D {\bf 43}, 4111-4120 (1991)
\item"\l HP\r"
Haagerup, U., Pisier, G.: Bounded linear operators between $C^*$-algebras.
Preprint 1993, to appear in Duke Math. J.
\item"\l Jim\r"
Jimbo, M.: Yang-Baxter equation in integrable systems. Advanced Series
in Mathematical Physics 10, Singapore: World Scientific 1989
\item"\l JSW1\r"
J\o rgensen, P.E.T., Schmitt, L.M., Werner, R.F.: $q$-Canonical
Commutation Relations and Stability of the Cuntz Algebra. Preprint 1992,
to appear in Pacific J. Math.
\item"\l JSW2\r"
J\o rgensen, P.E.T., Schmitt, L.M., Werner, R.F.:
Positive representations of general Wick ordering commutation relations.
Preprint 1993.
\item"\l Man\r"
Manin, Y.I.: Topics in Noncommutative Geometry.
Princeton: Princeton University
Press 1991
\item"\l Pau\r"
Paulsen, V.: Completely bounded maps and dilations. Pitna
Research Notes 146, New York: John
Wiley \& Sons 1986
\item"\l Pis1\r"
Pisier, G.: Multipliers and lacunary sets in non-amenable groups. Preprint
1992
\item"\l Pis2\r"
Pisier, G.: The Operator Hilbert Space OH, Complex Interpolation and Tensor
Norms. Preprint 1993, to appear in Memoirs AMS
\item"\l Sol\r"
Solomon, L.: The orders of the finite Chevalley groups. J. Algebra {\bf 3},
376-393 (1966)
\item"\l Spe1\r"
Speicher, R.: A New Example of \lq Independence' and \lq White Noise'.
Probab. Th. Rel. Fields {\bf 84}, 141-159 (1990)
\item"\l Spe2\r"
Speicher, R.: Generalized Statistics of Macroscopic Fields.
Lett. Math. Phys. {\bf 27}, 97-104 (1993)
\item"\l Voi\r"
Voiculescu, D.: Symmetries of some reduced free product
$C^*$-algebras. In:
Operator Algebras and their Connection with
Topology and Ergodic Theory
(LNM 1132, pp. 556-588), Heidelberg: Springer 1985
\item"\l VDN\r"
Voiculescu, D., Dykema, K., Nica, A.: Free Random Variables.
AMS 1992
\item"\l Was\r"
Wassermann, S.: Injective $W^*$-algebras. Math. Proc. Camb. Phil.
Soc. {\bf 82}, 39-47 (1977)
\item"\l Wen\r"
Wenzl, H.: Representations of braid groups and the quantum Yang-Baxter
equation. Pacific J. Math. {\bf 145}, 153-180 (1990)
\item"\l Zag\r"
Zagier, D.: Realizability of a Model in Infinite Statistics. Comm. Math.
Phys. {\bf 147}, 199-210 (1992)
\endroster
\enddocument